# Distributed Optimal Dynamic State Estimation for Cyber Intrusion Detection in Networked DC Microgrids

Tuyen V. Vu
ECE Department
Clarkson University
Potsdam, NY, USA
tvu@clarkson.edu

Bang Le
ECE Department
Clarkson University
Potsdam, NY, USA
nguyenbl@clarkson.edu

*Abstract*—In this paper, we present a novel distributed state estimation approach in networked DC microgrids to detect the false data injection in the microgrid control network. Each microgrid monitored by a distributed state estimator will detect if there is manipulated data received from their neighboring microgrids for control purposes. A dynamic model supporting the dynamic state estimation will be constructed for the networked microgrids. The optimal distributed state estimation, which is robust to load disturbances but sensitive to false data injected from neighboring microgrids will be presented. To demonstrate the effectiveness of the proposed approach, we simulate a 12kV three-bus networked DC microgrids in MATLAB/Simulink. Residual information corresponding to the false data injected from neighbors validates the efficacy of the proposed approach in detecting compromised agents of neighboring microgrids.

## I. INTRODUCTION

### A. Literature Review

Power outages can shut down critical infrastructures such as hospitals, water treatment plants, military services, and other emergency services. Financial consequences are significant. For example, disruptions in the U.S. electric power systems are estimated to be $25-70 billion annually [1]. During outages, microgrids can provide resilient energy service to critical infrastructure. While most power outages are the result of extreme weather events, there is increasing concern about outages caused by cyber-attacks [2]. Microgrids consist of distributed energy resources (DER) that provide power to local load devices. Effective operations of microgrids require advanced measurement, communication, and control via distributed controllers, sensors, actuators, and field devices. The measurement, communication, and control devices can be connected internally and externally via local area network or wide area network. Therefore, as a cyber-physical system (CPS), microgrids are particularly vulnerable to cyber-attacks due to their distributed nature and their critical resilience function. An extreme weather-induced outage may take weeks to restore and may cause significant economic and personal hardship [3]. Therefore, methodologies that improve microgrids' situational awareness of cyber-attacks are important.

In CPS, cyber intrusions are classified differently using different terms such as bias injection attack, zero dynamics attack, denial of service (DoS) attacks, eavesdropping attack, replay attack, stealthy attack, covert attack, and dynamic false data injection attacks [4]–[7]. However, all these attacks still focus on one or more components of CPS Data Confidentiality Integrity and Availability (CIA) triad, defined in common information security practices [8]. Individual types of attack have specific characteristics that influence the CIA-triad in individual ways. For example, while DoS attacks mainly affect the data availability, other attacks like replay, stealthy, dynamic false data injection and covert attacks influence data confidentiality and integrity. Attackers can manipulate the system control and management via 1) remote access to control system LAN network with poorly configured firewalls or 2) infected field devices [9].

As a CPS, microgrids can encounter the same types of attacks. DoS attacks can cause multiple issues to microgrids; however, once the DoS event occurs, the system operator very likely recognizes that the system is under attack. A more severe cyber-attack found in power systems is stealthy false data injection (FDI), where attackers corrupt the measurement and/or control data. In power transmission systems, the popular method used to detect bad measurement data is the static state estimator (SSE) based on weighted least squared (WLS). However, SSE can be manipulated by attackers if the power network topology is known [10]–[12]. Therefore, the cyber-attacks stay undetected as the attack indicator (control/measurement residuals) of the method is kept under predefined detected level. SSE approaches are not only susceptible to advanced cyber-attack policies they may also not be applicable for microgrids as there are more dynamical interactions among loads, generation, and distribution devices in microgrids.

There are methods for detecting FDI in microgrids in literature [13]–[17]. Recent relevant literature for FDI in microgrid can be found as follows: [13] investigated the FDI detection for consensus control of DC microgrids utilizing unknown input observer; however, the microgrid network model is described as quasi-static. Therefore, dynamic interactions within the microgrid can be omitted as the inductive interactions among devices are ignored. [14] proposed using invariants as fixed boundaries of voltages and currents without depending on the system model to detect the anomalies in a distributed control system of a DC microgrid; however; in many cases, a deviated voltage information, which still stays within the safety limit can arbitrarily drive the system into instability.

*B. Research Contributions*

In this paper, we 1) propose a novel optimal state estimation technique for FDI method in networked DC microgrids, 2) employ a dynamic model of the microgrid instead of the normally-used quasi-static model of DC microgrids for the design of optimal estimator, and 3) the proposed method is a distributed algorithm.

The remaining of the paper is organized as follows: Section II elaborates the proposed distributed state estimation methodology for networked DC microgrids. Section III provides case studies, and analyzes and discusses the results. Section IV concludes the achievements of the paper.

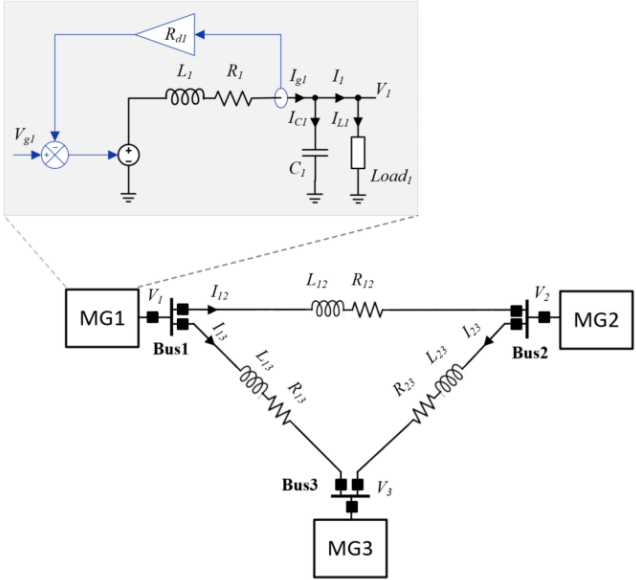

Figure 1. Representation of a three-bus DC microgrid circuit.

## II. DISTRIBUTED STATE ESTIMATION FOR NETWORK DC MICROGRIDS

In this section, the model of networked DC microgrids that support distributed state estimation will be elaborated. The model is then partitioned for our proposed distributed optimal state estimation.

*A. DC Microgrid Dynamic Modeling and Partitioning*

In this paper, we utilize a similar dynamic model developed in [18] for networked DC microgrids. The model we develop will include droop parameters as they are normally applied to control multi-terminal microgrids [19]. The dynamics of the system will reveal the nature of cyber-attacks/other physical incidents.

To have an equivalent model for networked DC microgrids, we will start with a three-bus networked DC microgrid as shown in Figure 1. Each microgrid can be modeled as an equivalent circuit with unknown internal load devices. In the figure, the equivalent circuit of microgrid 1 (MG1) is expanded for visualization. Equivalent droop control is implemented in each microgrid with the droop gain ($R_{d1}$). For simplification, dynamics of the equivalent source $V_{g1}$ is omitted in the modeling of the system. Based on that, MG1 can be modeled dynamically as follows:

$$\frac{dV_1}{dt} = \frac{1}{C_1}I_{g1} - \frac{1}{C_1}I_{12} - \frac{1}{C_1}I_{13} - \frac{1}{C_1}I_{L1}$$

$$\frac{dI_{g1}}{dt} = -\frac{1}{L_1}V_1 + \tau_1 I_{g1} + \frac{1}{L_1}V_{g1}$$

$$\frac{dI_{12}}{dt} = \frac{1}{L_{12}}V_1 + \tau_{12}I_{12} - \frac{1}{L_{12}}V_2$$

$$\frac{dI_{13}}{dt} = \frac{1}{L_{13}}V_1 + \tau_{13}I_{13} - \frac{1}{L_{13}}V_3, \tag{1}$$

where

$$\tau_1 = -\frac{R_{d1}+R_1}{L_1}, \tau_{12} = -\frac{R_{12}}{L_{12}}, \tau_{13} = -\frac{R_{13}}{L_{13}}.$$

(1) is the representation for MG1 coupling with MG2 and MG3. However, every microgrid can be generalized as a microgrid $i$, which interconnects with its neighbor $j$ ($j = 1, 2, \ldots, N_j$). The dynamic model of the a microgrid $i$ in the network is shown as

$$\dot{x}_i = A_{ci}x_i + B_{ci}u_i + E_{ci}d_i + \sum_{j=1}^{N_i} A_{cij}x_{ij} + w_i$$

$$y_i = C_{ci}x_i + v_i, \tag{2}$$

where $x_i$ is the local state vector with state matrix $A_{ci}$, $u_i$ is the control input vector with control matrix $B_{ci}$, $d_i$ is the disturbance with known disturbance structure $E_{ci}$, $w_i$ is the zero-mean process noise with covariance matrix $Q_i$, $y_i$ is the measurement output, and $v_i$ is the zero-mean measurement noise with covariance matrix $R_i$, $x_{ij}$ is the coupling states of neighbor $j$ with coupled-state matrix $A_{cij}$.

Take the three-bus networked microgrid as an example, the local state for MG1, $x_1 = [V_1, I_{g1}, I_{12}, I_{13}]^T$. The control variable, $u_1 = V_{g1}$, the disturbance, $d_1 = I_{L1}$, the measurement $y_1 = [V_1, I_{g1}, I_{12}, I_{13}]^T$, and the neighboring states $x_{1,j=1} = V_2$, $x_{1,j=2} = V_3$. The representation of the matrix components is detailed in the Appendix.

To detect the false data injected from neighbors, the neighboring states $x_{ij}$ will be referred to as additional inputs of system (2). Therefore, (2) can be rewritten as

$$\dot{x}_i = A_{ci}x_i + B_{cxi}u_{xi} + E_{ci}d_i + w_i$$
$$y_i = C_{ci}x_i + v_i, \qquad (3)$$

where $B_{cxi} = [B_{ci}\ A_{cij}]$, $u_{xi} = [u_i^T, (x_{ij})^T]^T$.
In the discrete domain, the system state is expressed as

$$x_{i,k+1} = A_i x_{i,k} + B_{xi} u_{xi,k} + E_i d_{i,k} + w_{i,k}$$
$$y_{i,k} = C_i x_{i,k} + v_{i,k}. \qquad (4)$$

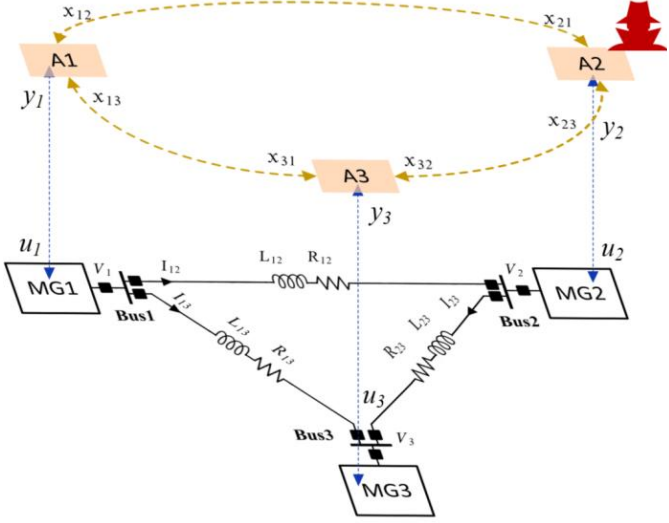

Figure 2. Distributed state estimation architecture in cyber-physical systems (networked microgrid), where a control agent is compromised.

This partitioned state-space system architecture will be utilized for the optimal state estimator proposed in the next sub-section.

### B. Distributed Optimal State Estimation

The proposed distributed optimal estimator for the three-bus example is shown in Figure 2. Each agent (A1, A2, or A3) exchange information $x_{ij}$ (bus voltage) with each other for the bus-voltage control purpose (Figure 3). Any agent or communication channel can be the target of cyber-attack. For example, A2 is compromised and start sending the false information to A1 and A3 to in order to disrupt the cyber-physical system.

An optimal state estimation that rejects the effects of unknown disturbances for detecting faults in a single control system was proposed in [20]. However, the method has not been analyzed for cyber-attack detection. In this paper, a distributed optimal estimation method is proposed for distributed FDI detection method. In this case, each distributed agent $A_i$ ($i = 1,2,3$) performs a distributed optimal state estimation (5) with information exchanged with their neighbors.

$$z_{i,k+1} = F_{i,k+1} z_{i,k} + T_i B_{xi}(u_{xi,k} + f_{uij,k}) + K_{i,k+1} y_{i,k}$$
$$\hat{x}_{i,k+1} = z_{i,k+1} + H_i y_{i,k+1}, \qquad (5)$$

where $a_{uij} = [0\ a_{ij}^T]^T$ in which $a_{ij}$ is the attack vector generated by neighbor $j$. With this observer, the state estimation error $e_{i,k+1} = x_{i,k+1} - \hat{x}_{i,k+1}$ in agent $i$ can be expressed as

$$e_{i,k+1} = F_{i,k+1} e_{i,k} - K_{i,k+1}^1 v_k - H_i v_{i,k+1}$$
$$+ (I - H_i C_i) w_i - (F_{i,k+1} - (I - H_i C_i)A_i + K_{i,k+1}^1 C_i) x_{i,k}$$
$$+ (I - H_i C_i) E_i d_{i,k} - [K_{i,k+1}^2 - F_{i,k+1} H_i] y_k$$
$$- [T_i - (I - H_i C_i)] B_{xi} u_{xi,k} - T_i B_{xi} a_{uij,k}, \qquad (6)$$

where $K_{i,k} = K_{i,k}^1 + K_{i,k+1}^2$.

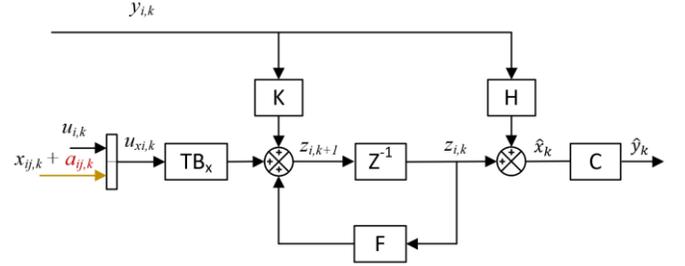

Figure 3. Optimal state estimation for distributed estimator $i$.

Define the following conditions:

$$(I - H_i C_i) E_i = 0$$
$$T_i - (I - H_i C_i) = 0$$
$$F_{i,k+1} - (I - H_i C_i) A_i + K_{i,k+1}^1 C_i = 0$$
$$K_{i,k+1}^2 - F_{i,k+1} H_i = 0, \qquad (7)$$

Solutions of (7) if exist will be

$$H_i = E_i [(C_i E_i)^T C_i E_i]^{-1} (C_i E_i)^T \qquad (8a)$$
$$T_i = (I - H_i C_i) \qquad (8b)$$
$$F_{i,k+1} = (I - H_i C_i) A_i - K_{i,k+1}^1 C_i \qquad (8c)$$
$$K_{i,k+1}^2 = F_{i,k+1} H_i, \qquad (8d)$$

where $K_{i,k+1}^1$ needs to be found for the stable and optimal observer. If these conditions hold, (6) will become

$$e_{i,k+1} = F_{i,k+1}e_{i,k} - K^1_{i,k+1}v_{i,k} - H_i v_{i,k+1}$$

$$+T_i w_i - T_i B_{xi} a_{uij,k}. \quad (9)$$

Based on this relationship, define the variance of $e_{i,k+1}$ as $P_{i,k+1}$, where

$$P_{i,k+1} = F_{i,k+1}P_{i,k}F^T_{i,k+1} + K^1_{i,k+1}R_{i,k}K^{1T}_{i,k+1}$$

$$-H_i R_{i,k+1}H_i^T + T_i Q_{i,k+1}T_i^T - T_i B_{xi}S_{uij,k}(T_i B_{xi})^T. \quad (10)$$

The parameter $S_{uij,k}$ is the unknown variance of $a_{uij,k}$. Therefore, without the attack, the variance will become

$$P_{i,k+1} = F_{i,k+1}P_{i,k}F^T_{i,k+1} + K^1_{i,k+1}R_{i,k}K^{1T}_{i,k+1}$$

$$-H_i R_{i,k+1}H_i^T + T_i Q_{i,k+1}T_i^T. \quad (11)$$

Therefore, the optimal observer can be achieved via minimization of $Trace(P_{i,k+1})$ with respect to $K^{1T}_{i,k+1}$. $\min_{K^1_{i,k+1}}(Trace((P_{i,k+1})))$ the estimator will be independent of disturbance, least dependent on measurement noise, and sensitive to the attack. $K^1_{i,k+1}$ can be found as

$$K^1_{i,k+1} = (I - H_i C_i)A_i P_{i,k}C_i^T (C_i P_{i,k}C_i^T + R_{i,k})^{-1}. \quad (12)$$

Therefore, the optimal observer parameters are iteratively obtained via (8), (11), and (12).

From (9), the residual $r_{i,k} = y_{i,k} - \hat{y}_{i,k}$ generated from the optimal observer will be

$$r_{i,k} = C_i e_{i,k+1} + v_{i,k}$$

$$r_{i,k} = C_i F_{i,k+1}e_{i,k} + (I - C_i K^1_{i,k+1})v_{i,k} - C_i H_i v_{i,k+1}$$

$$+C_i T_i w_i - C_i T_i B_{xi} a_{uij,k}. \quad (13)$$

Therefore, $r_{i,k}$ is most sensitive to the attack vector $a_{uij,k}$ and can be used to detect the attack event. In the next section, case studies will validate the feature of the proposed algorithm.

### III. CASE STUDIES

To demonstrate the proposed method, we utilize a 3-bus 12-kV DC system model. Each bus contains a microgrid, which has identical parameters shown in Table I. Each microgrid is rated at 50 MW and has an equivalent droop parameter of 5%. Process and measurement noise variances are shown in the Appendix.

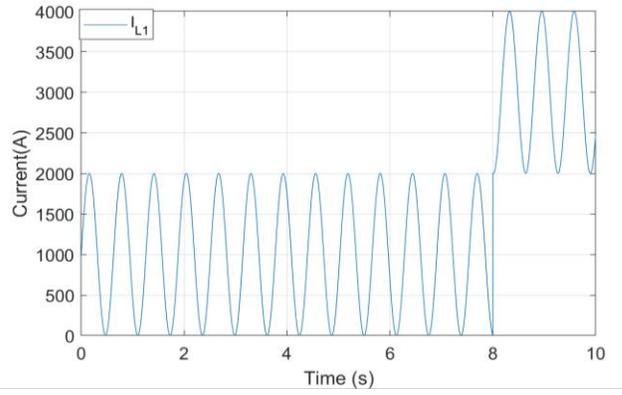

Figure 4. load profile of one MG.

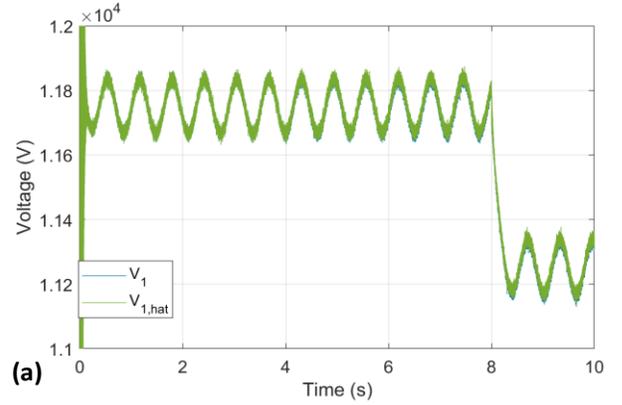

(a)

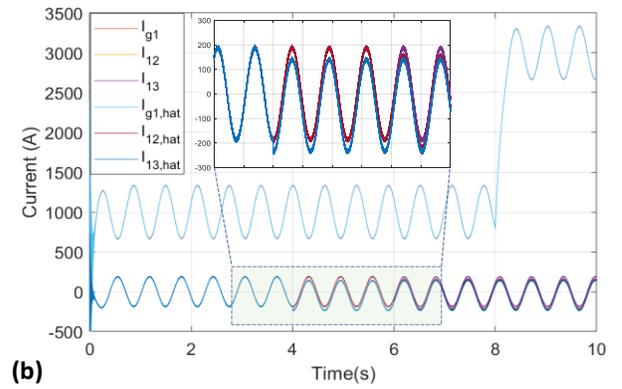

(b)

Figure 5. Distributed estimator $I$'s performance against load disturbance.

Table I. 12-kV DC Microgrid $i$'s parameters with interconnecting parameters

| Symbol | Description | Value |
|---|---|---|
| $R_i$ | Equivalent internal resistance | 0.05 Ω |
| $L_i$ | Equivalent internal inductance | 3 mH |
| $C_i$ | Equivalent output capacitance | 10 μF |
| $R_{ij}$ | Line resistance | 0.1 Ω |
| $L_{ij}$ | Line inductance | 0.5 mH |

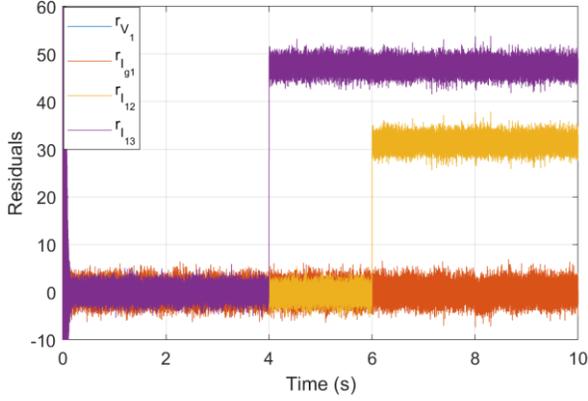

Figure 6. Residual generation.

Based on the parameters, we perform the following simulation case. We assume that a cyber attacker compromises distributed controller of MG2 and MG3 then sends false voltage data to MG1 at 6s and 4s with the corresponding biases of 100 V and 150 V, respectively. In this case, the load of each MG has a large range of current variation as shown in Figure 4. Estimated current and voltage of MG1 are shown in Figure 5. Residuals are shown in Figure 6. At 8s, the loads of MGs increase their current consumption by 2000 A.

As seen in Figure 5, the optimal state estimator eliminates the effect of load disturbance as prior to the attack events at 4s, the estimated values track the measured value. The residual indicators in Figure 6 also show the zero-mean residuals before 4s.

The first attack conducted by MG3 at 4s is detected as the estimated current (Figure 5b) $\hat{I}_{13}$ start to deviate from the measured current $I_{13}$. This is also shown in the residual component $r_{I_{13}} = I_{13} - \hat{I}_{13}$ as it changes to a non-zero component at 4s. As the current $I_{13}$ represents the coupling between MG1 and MG3, the MG1 could infer that the attack was conducted by MG3. Similarly, Figure 5b also shows that the estimated $\hat{I}_{12}$ starts deviating from the measured current $I_{12}$ at 6s. On the other hand, the residual $r_{I_{12}} = I_{12} - \hat{I}_{12}$ representing the coupling between MG1 and MG2 changes to a non-zero residual. This proves the attack at 6s originates from MG2.

At 8s, the load changes cause the increased consumption current and made the bus voltage drop; however, the residuals remained constant. This prove that the load disturbance at 8s is totally decoupled from the cyber-attacks.

Therefore, the proposed method effectively detects the attack in networked DC microgrids and is robust to load disturbance.

## IV. CONCLUSION

In this paper, we proposed a novel distributed state estimation methodology to detect cyber-attacks in distributed networked microgrids. A dynamic model was utilized for attack detection. Each MG will have the capability to detect compromised neighbor so that the control system can be informed for the system reconfiguration and attacks isolation. We demonstrated the proposed method in a three-bus microgrid network.

## APPENDIX

The model of MG1 and its coupling part is as:

$$A_1 = \begin{bmatrix} 0 & \frac{1}{C_1} & -\frac{1}{C_1} & -\frac{1}{C_1} \\ -\frac{1}{L_1} & \tau_1 & 0 & 0 \\ \frac{1}{L_{12}} & 0 & \tau_{12} & 0 \\ \frac{1}{L_{13}} & 0 & 0 & \tau_{13} \end{bmatrix}, B_1 = \begin{bmatrix} 0 \\ \frac{1}{L_1} \\ 0 \\ 0 \end{bmatrix}, E_i = \begin{bmatrix} -\frac{1}{C_1} \\ 0 \\ 0 \\ 0 \end{bmatrix}$$

$$C_1 = \begin{bmatrix} 1 & 0 & 0 & 0 \\ 0 & 1 & 0 & 0 \\ 0 & 0 & 1 & 0 \\ 0 & 0 & 0 & 1 \end{bmatrix}, A_{1,j=2} = \begin{bmatrix} 0 \\ 0 \\ 0 \\ -\frac{1}{L_{13}} \end{bmatrix}, A_{1,j=1} = \begin{bmatrix} 0 \\ 0 \\ -\frac{1}{L_{12}} \\ 0 \end{bmatrix}$$

$$Q_1 = diag\{q1, q2, q3, q4\}, R_1 = diag\{r1, r2, r3, r4\}, \quad (14)$$

where $n_1$ is the size of the state matrix $A_1$, $m_1$ is the number of measurements, and $q1$-$q4$ and $r1$-$r4$ are the variances of the process and measurement noise. Using the parameters indicated in Table I, parameters in (14) can be calculated except the variances. The process noise and measurement noise variances are assumed and calculated as

$$Q_1 = diag\{10,10,10,10\}, R_1 = diag\{100,100,10,10\}.$$